\documentclass{article}%
\usepackage{amsmath}
\usepackage{amsfonts}
\usepackage{amssymb}
\usepackage{graphicx}%
\providecommand{\U}[1]{\protect\rule{.1in}{.1in}}

\begin{document}

\title{On Resolving Problems with Conditionality and Its Implications for
Characterizing Statistical Evidence}
\author{Michael Evans$^{\ast}$ and Constantine Frangakis**\\Dept. of Statistics, University of Toronto$^{\ast}$ and \\Dept. of Biostatistics, Johns Hopkins University$^{\ast\ast}$}
\date{}
\maketitle

\begin{center}
\textbf{Abstract}
\end{center}

The conditionality principle $C$ plays a key role in attempts to characterize
the concept of statistical evidence. The standard version of $C$ considers a
model and a derived conditional model, formed by conditioning on an ancillary
statistic for the model, together with the data, to be equivalent with respect
to their statistical evidence content. This equivalence is considered to hold
for any ancillary statistic for the model but creates two problems. First,
there can be more than one maximal ancillary in a given context and this leads
to $C$ not being an equivalence relation and, as such, calls into question
whether $C$ is a proper characterization of statistical evidence. Second, a
statistic $A$ can change from ancillary to informative (in its marginal
distribution) when another ancillary $B$ changes, from having one known
distribution $P_{B},$ to having another known distribution $Q_{B}.$ This means
that the stability of ancillarity differs across ancillary statistics and
raises the issue of when a statistic can be said to be truly ancillary. It is
therefore natural, and practically important, to limit conditioning to the set
of ancillaries whose distribution is irrelevant to the ancillary status of any
other ancillary statistic. This results in a family of ancillaries for which
there is a unique maximal member. This also gives a new principle for
inference, the stable conditionality principle, that satisfies the criteria
required for any principle whose aim is to characterize statistical
evidence.\medskip

\noindent\textbf{Keywords and phrases}: ancillaries, stable ancillaries,
strong ancillaries, the laminal ancillary, stable conditionality principle

\section{Introduction}

The conditionality principle $C$ has played a puzzling role in attempts to
develop a frequentist theory of statistical inference. On the one hand it
seems intuitively obvious, and even a necessary component of such a theory.
But it also produces a significant ambiguity due to nonequivalent applications
for which there seems to be no easy solution in terms of determining which is
correct or even if any are correct. Attempts to ignore this problem, typically
by considering certain applications as equivalent, produces the somewhat
strange phenomenon that $C$, a frequentist principle, can lead to the
likelihood principle $L$ which precludes any frequentist inferences, see
Evans, Fraser and Monette (1986) and Evans (2013) for discussion of this.

The fact that $C$ is not an equivalence relation, which any valid
characterization of statistical evidence must be, calls into question the
justification for $C.$ This can be considered as a logical inconsistency in
the definition of $C$. Moreover, as will be shown, the ancillary status of a
statistic can change to being informative when the distribution of another
ancillary statistic changes. This raises the issue of whether the distribution
of such a statistic is truly irrelevant for inference, which can be considered
as a statistical inconsistency in the definition of $C.$

The purpose of this paper is to propose a resolution to these problems. It is
argued that a correct characterization of the ancillary concept requires the
restriction of the set of possible ancillaries for use to a subset and this is
based upon very natural statistical criteria. Once the restriction is made,
there is a unique maximal member of this subset and this becomes the ancillary
to use as it makes the maximal reduction in the set of possible data values to
compare the observed data to in the conditional model. We show that natural
statistical criteria lead to the set being the minimal ancillaries, whose
maximum is the laminal ancillary as labelled by the taxonomy of Basu (1959).

One could argue that this isn't much of an advance, particularly because the
laminal ancillary is often trivial, but we would counterargue that it is
significant because it shows that the other ancillaries, besides the laminal,
are \textit{ineligible} to be used in the conditioning step. This establishes
the validity of some form of $C$ for inference and this has broad
implications. In particular, the idea that $C$ together with the sufficiency
principle $S$ can lead to $L$, as discussed, for example, in Birnbaum (1962),
Evans, Fraser and Monette (1986), Evans (2013) and many others, is completely
avoided and this applies similarly to the argument that $C$ alone can produce
$L.$ Additionally, it leads to a new and uncontroversial principle that
combines $S$ and a modified $C$ that still permits frequentist considerations
for inferences.

In Section 2 the conditionality principle is discussed. In Section 3 we
introduce the statistical criterion that assesses whether an ancillary
statistic is unstable (can become informative) if one merely changes the
distribution of another ancillary statistic. We show how this connects to the
minimal and laminal ancillaries. In Section 4 a principle is introduced which
satisfies both $S$ and the new conditionality principle. This principle forms
an equivalence relation in the class of all inference bases and so is indeed a
valid partial characterization of statistical evidence, which was Birnbaum's
intention. The proofs of all propositions are placed in the Appendix.

The conditionality principle has attracted many authors some of whom have
attempted resolutions. The papers Basu (1959), Basu (1962), Cox (1958), Cox
(1971), Kalbfleisch (1975), Buehler (1982), Stigler (2001) and Ghosh, Reid and
Fraser (2010) all represent interesting contributions and there are many more
which can be found in the references of these papers. To the best of our
knowledge nobody has presented a forceful argument for the laminal ancillary
as being the natural resolution and that is the outcome of the discussion in
Section 3.

\section{Principles and Ancillaries}

All of the principles $S,C$ and $L,$ are applied to inference bases. An
inference base $I=(M,x)$ is comprised of a statistical model
\[
M=(\mathcal{X},\mathcal{B},\{P_{\theta,X}:\theta\in\Theta\}),
\]
where $\mathcal{X}$ is a sample space containing all possible values for the
observed data $x$ of random object $X,\mathcal{B}$ is a $\sigma$-algebra on
$\mathcal{X}$ and $\{P_{\theta,X}:\theta\in\Theta\}$ is a collection of
probability measures defined on $\mathcal{B}$ indexed by model parameter
$\theta\in\Theta.$ For inference, the assumption is made that there is a true
value of $\theta,$ say $\theta_{true},$ such that, before it is observed,
$x\sim P_{\theta_{true},X}.$ The goal, once $x$ is observed, is to make
inference about which of the possible values of $\theta\in\Theta$ corresponds
to $\theta_{true}$ and these inferences are based somehow on the ingredients
$I=(M,x).$ More generally, our interest is in some marginal parameter
$\psi=\Psi(\theta)$ that has a real-world interpretation and it is desired to
know the value $\psi_{true}=\Psi(\theta_{true})$ and this requires dealing
with so-called nuisance parameters. This more general problem is ignored here
except to say that the concept of conditioning on an ancillary for the model
is still relevant for that context. Birnbaum (1962) considered the set of all
inference bases and for inference bases $I_{1}$ and $I_{2}$ with essentially
the same model parameter (or bijective relabellings thereof), indicated that
these inference bases contain the same statistical evidence about the true
value of the model parameter by writing $Ev(I_{1})=Ev(I_{2}).$

An \textit{ancillary statistic} for the model $M$ is a map $A:(\mathcal{X}%
,\mathcal{B)\rightarrow}(\mathcal{A},\mathcal{C)}$ such that the marginal
probability measure induced by $A$ satisfies $P_{\theta,A}=P_{A}$ for every
$\theta\in\Theta.$ In other words $A$ is ancillary when its marginal
distribution is independent of the model parameter and it is then claimed that
the observed value of $A(x)$ contains no information about $\theta_{true}.$
More than this, simple examples, like the two measuring instruments example in
Cox (1958), suggest that for frequentist inferences the initial model $M$ in
$I=(M,x)$ be replaced by $M_{|A(x)}=\{P_{\theta,X}(\cdot\,|\,A(x)):\theta
\in\Theta\},$ where $P_{\theta,X}(\cdot\,|\,A(x)$ is the conditional
probability measure for $X$ given the value $A(x).$ The principle $C$ then
states $Ev(M,x)=Ev(M_{|A(x)},x).$ An ancillary $A$ is a \textit{maximal
ancillary} if, whenever $A^{\prime}$ is another ancillary and there exists a
function $h$ such that $A=h(A^{\prime}),$ then $h$ is effectively a 1-1
function. So, the set of possible data values $\{z:A(z)=A(x)\}$ that is
conditioned on via $C,$ when the value $A(x)$ of a maximal ancillary is
observed, cannot be made smaller without losing ancillarity.

It is natural to make the greatest possible reduction in the set of possible
sample values we use for inference and so a possible full statement of $C$
would be to condition on a maximal ancillary. When there is a unique maximal
ancillary this is uncontroversial. As Example 1 shows, however, there can be
several maximal ancillaries. In such a case there is an ambiguity concerning
which maximal ancillary to use when applying $C$ as, for two maximal
ancillaries $A_{1}$ and $A_{2},$ inference bases $(M_{|A_{1}(x)},x)$ and
$(M_{|A_{2}(x)},x)$ can lead to quite different inferences, see Example 2. It
is shown in Evans (2013) that the lack of a unique maximal ancillary implies
that $C$ is not an equivalence relation on the set of all inference bases and
therefore, as currently stated, it is not a correct characterization of
statistical evidence. Also, it is shown there that, if $\bar{C}$ is the
smallest equivalence relation containing $C,$ then $\bar{C}=L.$ Similarly the
smallest equivalence relation containing $S\cup C$, which is also not an
equivalence relation, satisfies $\overline{S\cup C}=L$ and this is what the
proof of Birnbaum's theorem proves. So the lack of a unique maximal ancillary
leaves open the question of whether or not $C,$ or some modification, is
indeed a valid statistical principle that should be employed in statistical work.

Basu (1959) defined a \textit{minimal ancillary }as any ancillary which is a
function of every maximal ancillary and showed that there is a unique
ancillary in the class of minimal ancillaries, called the \textit{laminal
ancillary,} which is maximal in this class. The following example illustrates
these concepts.\smallskip

\noindent\textbf{Example 1.}

Suppose $M$ consists of two distributions as provided in the Table
\ref{table1} together with the likelihood ratio (LR). Actually it is a range
of examples as $\epsilon$ is any value satisfying $0<\epsilon<1/64.$ For each
such case the minimal sufficient statistic (mss) is the identity which is not
the case if $\epsilon=0$. This implies that all the ancillaries are functions
of the mss and this will prove important for our later discussion.
\begin{table}[tbp] \centering
$%
\begin{tabular}
[c]{|c|c|c|c|c|c|c|c|}\hline
$x$ & $1$ & $2$ & $3$ & $4$ & $5$ & $6$ & $7$\\\hline
$\theta=\theta_{1}$ & $\frac{1}{8}+\epsilon$ & $\frac{1}{8}-\epsilon$ &
$\frac{1}{8}+2\epsilon$ & $\frac{1}{8}-2\epsilon$ & $\frac{1}{14}$ & $\frac
{2}{14}$ & $\frac{4}{14}$\\\hline
$\theta=\theta_{2}$ & $\frac{1}{16}-\epsilon$ & $\frac{3}{16}+\epsilon$ &
$\frac{3}{16}+4\epsilon$ & $\frac{1}{16}-4\epsilon$ & $\frac{2}{14}$ &
$\frac{1}{14}$ & $\frac{4}{14}$\\\hline
$LR$ & $\frac{1/8+\epsilon}{1/16-\epsilon}$ & $\frac{1/8-\epsilon
}{3/16+\epsilon}$ & $\frac{1/8+2\epsilon}{3/16+4\epsilon}$ & $\frac
{1/8-2\epsilon}{1/16-4\epsilon}$ & $\frac{1}{2}$ & $2$ & $1$\\\hline
\end{tabular}
\ $%
\caption{Distributions in Example 1 together with likelihood ratios.}\label{table1}%
\end{table}%

Since any 1-1 function of an ancillary is ancillary, it is equivalent to
present all the preimage partitions induced by such statistics when
considering the ancillary structure of this model and some of these are
provided in the following table. It is clear from this table that the maximal
ancillaries are given by $A_{1}$ and $A_{2},$ as these give the finest
ancillary partitions, and so the laminal ancillary must be $L$ as it is the
finest partition containing both maximal ancillaries. The minimal ancillaries
are given by $\{T,B_{1},B_{2},B_{3},L\},$ where $T$ is the trivial ancillary,
as these are all coarsenings of both $A_{1}$ and $A_{2}$ and are presented in
Table \ref{table2}.%
\begin{table}[tbp] \centering
$%
\begin{tabular}
[c]{|c|c|}\hline
ancillary & partition of $\mathcal{X}$\\\hline
$T$ & $\{1,2,3,4,5,6,7\}$\\\hline
$B_{1}$ & $\{1,2,3,4,5,6\},\{7\}$\\\hline
$B_{2}$ & $\{1,2,3,4,7\},\{5,6\}$\\\hline
$B_{3}$ & $\{1,2,3,4\},\{5,6,7\}$\\\hline
$L$ & $\{1,2,3,4\},\{5,6\},\{7\}$\\\hline
$A_{1}$ & $\{1,2\},\{3,4\},\{5,6\},\{7\}$\\\hline
$A_{2}$ & $\{1,3\},\{2,4\},\{5,6\},\{7\}$\\\hline
\end{tabular}
\ $\caption{The minimal ancillaries in Example 1.}\label{table2}%
\end{table}
There are ancillaries that are coarsenings of single maximal ancillaries such
as%
\begin{align*}
C_{1}  &  :\{1,3\},\{2,4\},\{5,6,7\}\\
C_{2}  &  :\{1,3,5,6\},\{2,4\},\{7\}
\end{align*}
which are coarsenings of $A_{2}$ but not of $A_{1}$ and there are many others.

If the sample space were shrunk to $\{1,2,3,4\},$ with the $1/2$ probability
for $\{5,6,7\}$ redistributed equally among the 4 sample points, then the
laminal ancillary becomes the trivial ancillary and this is not uncommon, as
noted in Basu (1959) where conditions for this to occur are discussed.
$\blacksquare$\smallskip

The following example demonstrates the ambiguity that a nonunique maximal
ancillary can produce and is adapted from Evans (2015).\smallskip

\noindent\textbf{Example 2. }

Consider the model given by Table \ref{table3} and suppose $x=1$ is observed.
The MLE\ of $\theta$ is $\hat{\theta}(1)=\theta_{1}.$
\begin{table}[tbp] \centering
$%
\begin{tabular}
[c]{|c|c|c|c|c|}\hline
$x$ & $1$ & $2$ & $3$ & $4$\\\hline
$\theta=\theta_{1}$ & $1/6$ & $1/6$ & $2/6$ & $2/6$\\\hline
$\theta=\theta_{2}$ & $1/12$ & $3/12$ & $5/12$ & $3/12$\\\hline
\end{tabular}
\ \ $\caption{Distributions in Example 2.}\label{table3}%
\end{table}
There are two maximal ancillaries as given by their partitions, namely
$A_{1}=\{\{1,2\},\{3,4\}\}$ and $A_{2}=\{\{1,3\},\{2,4\}\}.$ The sampling
distributions of the MLE obtained by conditioning on the maximal ancillaries
are as displayed in Table \ref{table4}.%
\begin{table}[tbp] \centering
$%
\begin{tabular}
[c]{|c|c|c|}\hline
& $\theta=\theta_{1}$ & $\theta=\theta_{2}$\\\hline
$P_{a}(\hat{\theta}(X)=\theta\,|\,A_{1}=\{1,2\})$ & $1/2$ & $1/2$\\\hline
$P_{b}(\hat{\theta}(X)=\theta\,|\,A_{1}=\{1,2\})$ & $1/4$ & $3/4$\\\hline
$P_{a}(\hat{\theta}(X)=\theta\,|\,A_{2}=\{1,3\})$ & $1/3$ & $2/3$\\\hline
$P_{b}(\hat{\theta}(X)=\theta\,|\,A_{2}=\{1,3\})$ & $1/6$ & $5/6$\\\hline
\end{tabular}
$\caption{Conditional distributions of the MLE in Example 2.}\label{table4}%
\end{table}
As can be seen, these sampling distributions are quite different and it is not
clear which to use as part of quantifying the uncertainty in the estimate.
$\blacksquare$

\section{Stable and Strong Ancillaries}

Despite the rich structure of the ancillary statistics, standard evidence
theory assumes (through the standard conditionality principle $C$) that
conditioning on different ancillary statistics is equally valid. We challenge
this assumption through two main perspectives, which give rise to a
resolution.\smallskip

\textbf{Reproducing the structure with a single maximal ancillary} As noted in
Evans (2013), the fact that more than one maximal ancillary can exist results
in $C$ not forming an equivalence relation on the set of all inference bases.
If we want to claim that a given principle does properly characterize when two
inference bases contain the same amount of statistical evidence concerning an
unknown $\theta,$ then it seems clear that the principle must induce an
equivalence relation. Therefore, $C$ needs to be modified if it is desirable
for conditioning on ancillaries to play a role in inference.

Basu (1959) introduced the concept that two ancillary subsets $A,B\in
\mathcal{B}$ for model $M$ \textit{conform} when $A\cap B$ is also ancillary.
The set of all ancillary subsets that conform to every other ancillary subset
is denoted by $\Gamma_{0}$ and it is proved that $\Gamma_{0}$ is a $\sigma
$-algebra and moreover this is the laminal ancillary $\sigma$-algebra in the
sense that it is the largest $\sigma$-algebra contained in all the $\sigma
$-algebras induced by the individual maximal ancillaries. This is effectively
saying that (allowing for 1-1 equivalences) the laminal ancillary statistic is
a function of every maximal ancillary. A further implication of this is that
the laminal ancillary $\sigma$-algebra is the largest minimal ancillary
$\sigma$-algebra and so the laminal ancillary statistic is the maximal minimal
ancillary statistic. Also, if there is a unique maximal ancillary then this is
also the laminal ancillary. This points to a special role for the laminal
ancillary especially since the laminal ancillary always exists and a
conditionality principle that prescribed conditioning on the laminal forms an
equivalence relation on the set of inference bases, see Section 4.

Although logical, this role has not been explored. Perhaps this is because the
laminal doesn't often produce a meaningful reduction. But also Basu's
development, while logical, doesn't provide a good statistical reason to adopt
the laminal as the logical ancillary to condition on. It is argued here,
however, that there is a key element that can be added to the story and with
this addition the laminal is not only a logical resolution, but is a
statistical necessity.\smallskip

\textbf{Addressing the transition of ancillaries to informative statistics
}The key idea in this development is the supposed irrelevance of the
distribution of an ancillary that is to be conditioned on. For after all, as
far as inference goes, this distribution plays absolutely no role whatsoever.
The statistical intuition behind this is that the distribution of the
ancillary is free of the parameter and so an observation from it contains no
information about $\theta_{true}.$ As such, it must be the case that, no
matter what distribution is assumed for an ancillary this cannot change the
basic information structure of the problem. Note that this is a more severe
requirement for what it means for a statistic to be ancillary. Two definitions
that capture this idea are now provided and their equivalence proved. It is
then proved that the set of ancillaries which satisfy this criterion has a
maximal member and it is the laminal ancillary. To avoid a measure-theoretic
presentation via $\sigma$-algebras, as in Basu (1959), it will be assumed here
that all ancillaries are discretely distributed on $%
\mathbb{N}
$ and that there are at most a countable number of ancillaries, as this is
sufficient for conveying the key ideas.

For ancillary $U$ for model $M,$ the following notation is adopted
\[
M=\sum_{i}P_{U}(\{i\})M_{|U=i}.
\]
This expresses the idea that the model $M$ is a mixture of the component
models obtained by conditioning on $U=i$ where the mixture probabilities are
given by the marginal distribution of $U.$ The following definitions capture
the idea that the distribution of $U$ should be irrelevant for the inference
problem. $\smallskip$

\noindent\textbf{Definition} An ancillary $U$ for model $M$ is called a
\textit{stable ancillary} for $M$ if, whenever $V$ is ancillary for $M,$ then
$U$ is ancillary for the mixture $\sum_{i}p_{i}M_{_{|V=i}}$ for every
probability distribution $(p_{1},p_{2},\ldots)$ on the set of possible values
for $V.$ An ancillary $U$ for model $M$ is called a \textit{strong ancillary}
for $M$ if any ancillary $V$ for $M$ is also ancillary for the mixture
$\sum_{i}p_{i}M_{_{|U=i}}$ for every probability distribution $(p_{1}%
,p_{2},\ldots)$ on the set of possible values for $U.\smallskip$

\noindent So $U$ is a stable ancillary when changing the distribution of any
other ancillary has no effect on the ancillarity of $U$ and $U$ is a strong
ancillary if changing the distribution of $U$ has no effect on the ancillarity
of any other ancillary. For any ancillary $U$ that is not stable, then
conditioning on the value of some other ancillary renders the value $U(x)$
informative which contradicts the underlying motivation that the value of an
ancillary statistic contains no evidence concerning $\theta_{true}$.
Similarly, if $U$ is not strong, then conditioning on the value $U(x)$ renders
the value of some other ancillary informative. Accordingly, it is difficult to
accept the claim that the value of an ancillary that is not stable/strong is
noninformative with respect to $\theta_{true}.$

In actuality, a stable ancillary is strong and a strong ancillary is stable as
the following result shows.$\smallskip$

\noindent\textbf{Proposition 1.} $U$ is a strong ancillary for $M$ iff it is a
stable ancillary for $M.\smallskip$

\noindent Given that stable and strong ancillaries are\ just different
expressions of the same concept, these will be referred to hereafter as stable ancillaries.

In part (i) of the following result it is now shown that a stable ancillary is
a minimal ancillary and a minimal ancillary is a stable ancillary. Since Basu
(1959) proved that the laminal ancillary is the maximal minimal ancillary this
establishes that the laminal ancillary is the maximal stable ancillary and,
for the sake of completeness, this is proved in part (ii).$\smallskip$

\noindent\textbf{Proposition 2.} (i) A stable/strong ancillary is a minimal
ancillary and conversely. (ii) There exists a maximal minimal ancillary (the
laminal ancillary).$\smallskip$

\noindent Since the word minimal doesn't really convey the positive aspects of
such ancillaries these will be referenced as stable ancillaries hereafter.

It is worth noting that the structure given by the minimal and laminal
ancillaries is really the largest ancillary structure within the model that
replicates the situation where there is a single maximal ancillary and, as
such, there is no ambiguity about which ancillary to condition on. This
coherence points to the laminal ancillary as playing a special role and this
is reinforced by the notion of stability of an ancillary.

The following example demonstrate numerically the extent to which, having an
incorrect distribution of an unstable ancillary (i) can transform another
unstable ancillary to informative; yet (ii) preserves the ancillary state of a
stable ancillary.$\smallskip$

\noindent\textbf{Example 3.}

Consider again Example 1 with $\epsilon=0.01$, but now consider what happens
to the ancillary state of the unstable ancillary $C_{2}$ and the stable
ancillary $L$, when the distribution of the unstable ancillary $A_{1}$ is
changed from $P_{A_{1}},$ as given by $(1/4,1/4,3/14,4/14),$ to a true
distribution that is unknown to the researcher, $P_{A_{1}}^{unknown},$ as
given by $(7/100,13/100,27/100,53/100)$, see Figure 1. It is then observed
that $L$ stays ancillary, as theory assures, namely, for both $\theta
=\theta_{1}$ and $\theta=\theta_{2},$ the distribution of $L$ is
$(1/2,3/14,4/14)$ under the first scenario and $(20/100,27/100,53/100)$ under
the second. However, the likelihood ratios of $C_{2},$ $P_{\theta_{1}}(C_{2}=$
a given value$)/P_{\theta_{2}}(C_{2}=$ a given value$),$ are largely away from
$1;$ $C_{2}$ has lost its ancillary state and is now informative.%
\begin{figure}[ptb]%
\centering
\includegraphics[
height=2.2814in,
width=4.4166in
]%
{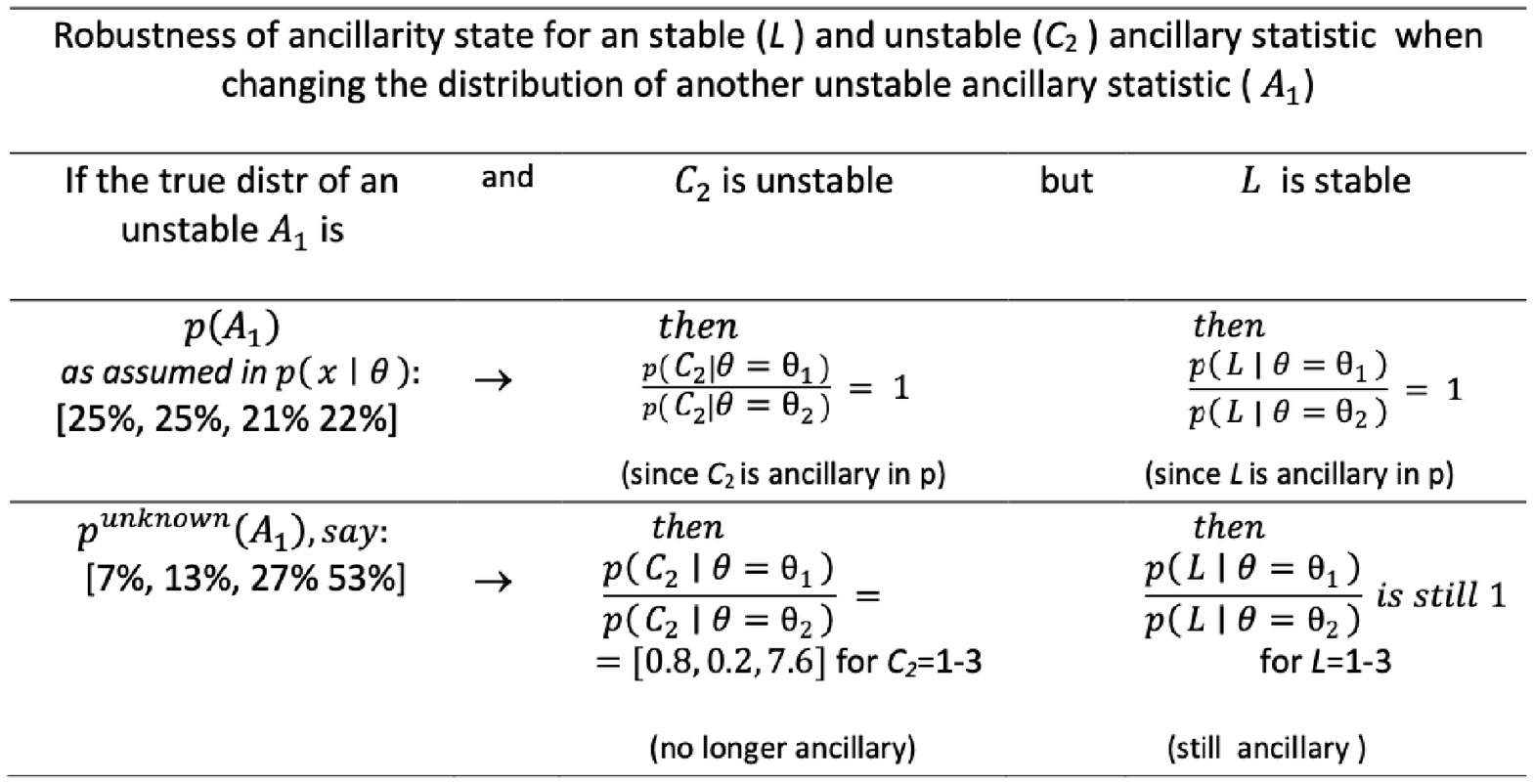}%
\caption{The result of changing the distribution of a nonstable ancillary in
Example 3.}%
\end{figure}

One may consider reasonable that such sensitivity of the ancillary state for a
statistic suggests that its ancillarity is not a structural feature of the
design, but is rather an erroneous coincidence. This possibility, while not
testable within the model, suggests that one should focus any conditioning
only on stable ancillaries. $\blacksquare\smallskip$

To see additionally why $C$ needs to be modified we examine the motivation for
conditioning as part of the inference process. This arises from considering
mixture experiments. Suppose there are a set of models say $\{M_{a}%
:a\in\mathcal{A}\},$ with $M_{a}=(\mathcal{X},\{P_{\theta,a}:\theta\in
\Theta\}),$ where the data $x$ will arise from one of these models. The model
that produces the data is obtained via a randomization procedure where a value
$a$ is produced with probabilities given by $P_{A}(\{a\})=P(A=a),$ on
$\mathcal{A}.$ This mixing produces the overall model $M=\sum_{a\in
\mathcal{A}}P_{A}(\{a\})M_{a}$ and $A$ is ancillary for $M.$ If the value of
$A=a_{0}$ is observed, then $C$ says that the inference base $(M_{a_{0}},x)$
is the one that is relevant for inference about $\theta.$ This seems
uncontroversial and therein lies the appeal of $C.$

The controversy surrounding $C$ arises when, rather than being presented with
a physical randomization device as part of a two-stage procedure, as just
described, we are presented with the inference base $(M,x)$ with $A$ being
ancillary for $M.$ Since $M$ can be at least be formally considered as a
mixture model via $A,$ it then seems reasonable to replace $(M,x)$ by
$(M_{a_{0}},x),$ where $A(x)=a_{0},$ for inference about $\theta.$

But now consider two studies conducted by statisticians 1 and 2 concerning the
true value of the quantity $\theta$ but suppose different randomization
schemes are used in each. So, in the $i$-th study the collection of models is
given by $\{M_{ia}:a\in\mathcal{A}_{i}\}$ and the relevant ancillary is
$A_{i}.$ Suppose that the results of the mixing produces the same overall
model $M$ and furthermore the same data $x$ is obtained. This may seem
unrealistic, but recall that in the end this is the situation that confronts
us when considering a model with multiple ancillaries and we wish to justify
conditioning on one of them.

It would seem then that both studies would conclude that the evidence about
the true value of $\theta$ in the inference base $(M,x)$ is the same but the
expression of this will be different, and result in different conditional
inference bases, unless effectively the same maximal ancillary is being used
for the mixing. In Example 1, suppose the two randomization schemes are
specified by the maximal ancillaries $A_{1}$ and $A_{2}$ as this will be a
case where the conditional inference bases will be different. Recall, however,
that the specific distributions for the $A_{i}$ are supposedly irrelevant for
inference about $\theta$ and indeed these play no role in the actual
inferences. But now suppose, for whatever reason, statistician 1 decides to
modify their randomization scheme by changing the distribution of $A_{1}$ say
from $P_{A_{1}}$ to $P_{A_{1}}^{\prime}.$ This does not change the submodels
$M_{1a}$ and so this change in the ancillary distribution seems innocuous to
statistician 1 as their inferences will not change due to the irrelevance of
the distribution of the ancillary. The overall model $M,$ however, has changed
to $M^{\prime}$ and this may produce a conflict with statistician 2 because it
may be that $A_{2}$ is no longer ancillary in $M^{\prime}$ and is now
informative. Statistician 2 can now rightly claim that the distribution of
$A_{1}$ is definitely relevant to the inference process and so there is a
contradiction between the two statisticians.

This demonstrates that there is a clear contradiction that resides within the
reasoning that justifies $C,$ at least as long as it is silent about which
ancillaries are appropriate for the conditioning step$.$ The content of this
paper has demonstrated how to resolve this contradiction by making sure that
any ancillaries that are used do not produce the phenomenon just described.
The relevant ancillaries to use are the stable ancillaries and indeed their
marginal distributions are irrelevant for inference. The irrelevance of the
marginal distribution of a stable ancillary is similar to the irrelevance of
the conditional distribution of the data given a mss and both can be discarded
for inference. This recovers conditioning on an ancillary as a valid part of
the inference process. Of course, we want to make the maximal reduction via
conditioning, to eliminate as much of the variation as possible that has
nothing to do with $\theta,$ and this leads to conditioning on the laminal.

\section{Stable Conditionality and Evidence}

In discussing statistical evidence Birnbaum (1962) introduced the $Ev$
function defined on the set of all inference bases. When two inference bases
$I_{1},I_{2}$ were considered to be equivalent with respect to their content
of statistical evidence, this was denoted by $Ev(I)=Ev(I_{2}).$ Birnbaum did
not, however, specify the value of $Ev(I).$ While this is understandable, this
approach is modified here as evidence functions are fully defined (up to 1-1
equivalence due to relabellings) for the principles discussed. The basic
reason for this is that a principle of inference should not only state an
equivalence, but also prevent the usage of aspects of an inference base that
are identified as irrelevant for the inference process. As pointed out in
Durbin (1976), ensuring that this didn't happen was one way of preventing
Birnbaum's proof of his well-known theorem. We still do not give a full
definition of $Ev$ but it is argued that this takes us some steps closer and
that such restrictions are a necessity.

In what follows, we examine the consequences that arise for statistical
evidence as described in Birnbaum, if one focuses on the set of stable
ancillaries that are functions of a mss for a model $M$, namely,%
\begin{equation}
\mathcal{A}_{M}=\{A:A\text{ is a stable ancillary and a function of a mss for
model }M\}. \label{focus}%
\end{equation}
It was pointed out in Durbin (1970) that restricting to ancillaries that are
functions of a mss voided the proof of Birnbaum's theorem. Evans et al. (1986)
argued that this was a natural restriction because otherwise the information
being conditioned on via the ancillary was precisely the information being
discarded as irrelevant via sufficiency in Birnbaum's proof. As such, there
existed a contradiction between the principles $S$ and $C$ in that context$.$
The restriction to ancillaries that are functions of a mss also seems implicit
in Fisher's development of the ancillarity concept, as documented in Stigler (2001).

Based on the developments in Section 3, the restriction is made to those
ancillaries that are stable\ because these are in a sense the ancillaries that
truly introduce no information into the analysis concerning the true
distribution. It is to be noted that there still is a place in a statistical
analysis for ancillaries that are not functions of a mss as, for example, in
regression analysis with normal error where the standardized residuals are
ancillaries that are not functions of the mss but play a key role in model
checking. Our concern here, however, is with the inference step and the
restriction to (\ref{focus}) seems essential in that context.

For simplicity, we suppose that the parameter space $\Theta=\{\theta
_{1},\theta_{2},...,\theta_{m}\}$ and the sample space $\mathcal{X}%
=\{x_{1},,...,x_{n}\}$ are both finite as this doesn't change the essential
meaning of the principles. Also we take $\mathcal{B}=2^{\mathcal{X}},$ the
power set of $\mathcal{X},$ and suppress this in the notation hereafter. It is
assumed that $\Theta$ is the same in any two inference bases that we consider
related via $Ev$ although it is possible to allow one parameter space to be a
1-1 relabelling of the other but this is ignored here.\ Also, it will always
be assumed that, for each $x_{i}\in\mathcal{X}$ then there is at least one
$\theta\in\Theta$ such that $P_{\theta}(\{x_{i}\})>0$ so the sample space
$\mathcal{X}$ cannot be made smaller.

A sufficient statistic $T$ is any function defined\ on $\mathcal{X}$ such
that, if $T(x)=T(y),$ then $x$ and $y$ are in the same equivalence class
associated with the sufficiency equivalence relation on $\mathcal{X}$ given
by\ $x\equiv_{S}y$ whenever there is a constant $c$ such that $P_{\theta
,X}(\{x\})=cP_{\theta,X}(\{y\})$ for every $\theta\in\Theta$. A mss is a
sufficient statistic $T$ such that when $x\equiv_{S}y,$ then $T(x)=T(y)$ and
so it is any function on $\mathcal{X}$ that indexes the equivalence classes.
The value of the mss represents the maximal reduction in the observed data
that results in no information loss concerning $\theta_{true}.$ A canonical
representative of the mss is, as discussed in Evans (2015), Lemma 3.3.2, given
by $T(x)=[x]$ where $[x]\subset\mathcal{X}$ is the equivalence class induced
by $\equiv_{S}$ on $\mathcal{X}.$ Any function on $\mathcal{X}$ that is
constant on each set $[x]$ and different on $[x]$ and $[y]$ when
$[x]\neq\lbrack y],$ can also serve as a mss. For example, when there is
$\theta_{i}\in\Theta$ such that $P_{\theta_{i},X}(\{x\})>0$ for all
$x\in\mathcal{X},$ then the mss can be taken to be
\[
T(x)=(P_{\theta_{1},X}(\{x\})/P_{\theta_{i},X}(\{x\}),\ldots,P_{\theta_{n}%
,X}(\{x\})/P_{\theta_{i},X}(\{x\})).
\]
Let $T:\mathcal{X}\overset{onto}{\mathcal{\rightarrow}}\mathcal{T}$ denote the
mss, however it is chosen, with model $M_{T}=(\mathcal{T},\{P_{\theta
,T}:\theta\in\Theta\}).$

The following statement of the sufficiency principle is equivalent to the
statement in Birnbaum (1962) but it is easier to use this version to prove
that $S$ is indeed an equivalence relation on the set of all inference bases,
see Evans (2015), Lemma 3.3.3. Here we allow for any version of the mss as
$h(T)$ where $h$ is a 1-1 function (a \textit{relabelling})\ defined on
$\mathcal{T}.$ This allows for relating two inference bases $(M_{1},x_{1})$
and $(M_{2},x_{2})$ that may have very different models but their minimal
sufficient statistics are essentially equivalent under such a relabelling and
so the principle is defined as a relation on the set of all inference
bases.\smallskip

\noindent\textbf{Sufficiency Principle }$(S)$ The\textbf{ }inference bases
$(M_{1},x_{1})$ and $(M_{2},x_{2}),$ with minimal sufficient statistics
$T_{1}$ and $T_{2}$ respectively, are equivalent under $S$ whenever there is a
a 1-1 onto, function $h:\mathcal{T}_{2}\mathcal{\rightarrow T}_{1}$ such that
$T_{1}=h\circ T_{2}$ and
\[
(M_{1,T_{1}},T_{1}(x_{1}))=(M_{2,h(T_{2})},h(T_{2}(x_{2}))).
\]
So when $(M_{1},x_{1})$ and $(M_{2},x_{2})$ are related via $S,$ the sampling
distributions of $T_{1}$ and $T_{2}$ are essentially the same as are the
observed values of these statistics. For example, as a particular application,
if model $(\mathcal{X},\{P_{X\mid\theta}:\theta\in\Theta\})$ has mss $T,$ then
observations $x,y\in\mathcal{X}$ satisfying $T(x)=T(y),$ together with the
model, contain the same evidence about $\theta_{true},$ i.e.,
\begin{equation}
Ev(\mathcal{X},\{P_{\theta,X}:\theta\in\Theta\},x)=Ev(\mathcal{X}%
,\{P_{\theta,X}:\theta\in\Theta\},y)\nonumber
\end{equation}
where the function $h$ is just the identity in this case.

While no image space is defined for $Ev$ it is necessary to do this for a
specific principle so that it is clear that the goal of the principle is to
also \textit{exclude} ingredients that are really extraneous to the intent of
the principle. It is immediate from $S$ that
\[
Ev(\mathcal{X},\{P_{\theta,X}:\theta\in\Theta\},x)=Ev(\mathcal{T}%
,\{P_{\theta,T}:\theta\in\Theta\},T(x))
\]
and this is undoubtedly the most important application of the principle,
namely, all inferences about the true value of $\theta$ are based on the model
for a mss and its observed value. This leads to the definition of the
\textit{minimal sufficiency evidence function} $Ev_{MS}$ given by
\[
Ev_{MS}(\mathcal{X},\{P_{\theta,X}:\theta\in\Theta\},x)=(\mathcal{T}%
,\{P_{\theta,T}:\theta\in\Theta\},T(x))=(M_{T},T(x)),
\]
for say the canonical mss $T,$ although any other equivalent version of the
mss could be used$.$ In other words, we are restricting what we consider an
appropriate presentation of the evidence based on $S.$ The ultimate evidence
function, whatever it may be, will be composed with $Ev_{MS}.$

For ancillary statistic $A$ for model $M=(\mathcal{X},\{P_{\theta,X}:\theta
\in\Theta\})$ we write $M_{|A(x)}=(\mathcal{X},\{P_{\theta,X\mid A(x)}%
:\theta\in\Theta\})$ for the family of derived conditional distributions on
$\mathcal{X}$ obtained by conditioning on the event specified by $A(x).$ The
discussion in Section 3 about ancillarity then leads to the following modified
conditionality principle where again we state a general version of the
principle that can be applied to relate (or not) any inference
bases.\smallskip

\noindent\textbf{Stable Conditionality Principle }$(SC)$\textbf{ }The\textbf{
}inference bases $(M_{1},x_{1})$ and $(M_{2},x_{2}),$ with minimal sufficient
statistics $T_{i}$ and laminal ancillaries $L_{i}\in\mathcal{A}_{M_{i}}$
respectively, are equivalent under $SC,$ whenever there is a a 1-1 onto,
function $h:\mathcal{T}_{2}\mathcal{\rightarrow T}_{1}$ such that
$T_{1}=h\circ T_{2}$ and
\begin{equation}
(M_{1,T_{1}\mid L_{1}(x_{1})},T_{1}(x_{1}))=(M_{2,h(T_{2})\mid L_{2}(x_{2}%
)},h(T_{2}(x_{2}))). \label{screl}%
\end{equation}
\smallskip\noindent For example, if model $(\mathcal{X},\{P_{X\mid\theta
}:\theta\in\Theta\})$ has mss $T$ and laminal ancillary $L\in\mathcal{A}_{M},$
then observations $x,y\in\mathcal{X}$ satisfying $T(x)=T(y),$ together with
the conditional model, contain the same evidence about $\theta_{true},$ i.e.,
\begin{equation}
Ev(\mathcal{X},\{P_{\theta,X|L(x)}:\theta\in\Theta\},x)=Ev(\mathcal{X}%
,\{P_{\theta,X|L(y)}:\theta\in\Theta\},y)\nonumber
\end{equation}
where the function $h$ is just the identity in this case.

It follows from $SC$ that
\[
Ev(\mathcal{X},\{P_{\theta,X}:\theta\in\Theta\},x)=Ev(\mathcal{T}%
,\{P_{\theta,T|L(x)}:\theta\in\Theta\},T(x))
\]
and this is undoubtedly the most important application of the principle. This
leads to the definition of the \textit{stable conditionality evidence
function} $Ev_{MS}$ given by
\begin{equation}
Ev_{SC}(\mathcal{X},\{P_{\theta,X}:\theta\in\Theta\},x)=(\mathcal{T}%
,\{P_{\theta,T|L(T(x))}:\theta\in\Theta\},T(x)) \label{stabev}%
\end{equation}
for say the canonical mss $T$ although any other equivalent version of the mss
could be used.

It is necessary to prove that $SC$ is an equivalence relation on the set of
all inference bases as part of establishing that $Ev_{SC}$ is a valid
characterization of statistical evidence.\smallskip

\noindent\textbf{Proposition 3.} $SC$ is an equivalence relation on the set of
inference bases.\smallskip

\noindent It is obvious that, as relations on the set of all inference bases,
$SC\subset C.$ The fact that $SC$ is an equivalence relation establishes that
this containment is proper because it has been established that $C$ is not an
equivalence relation, see Evans (2013) or Evans (2015), Lemma 3.3.4. It has
also been shown in these references that the smallest equivalence relation
containing $C$ is $L.$ So an interesting consequence of Proposition 3 is that
$L$ cannot be obtained from $SC$ in this way.

Similarly, the same references establish that the relation given by $S\cup C$
is not an equivalence relation and the proof of Birnbaum's Theorem establishes
that the smallest equivalence relation containing $S\cup C$ is $L.$ In this
case the following establishes that $S\subset SC$ so Birnabum's Theorem does
not follow from $S$ and $SC.$\smallskip

\noindent\textbf{Proposition 4.} As relations on the set of all inference
bases $S\subset SC.$\smallskip

\noindent Note that $SC$ only requires that the conditional models
$M_{i,T_{i}\mid L_{i}(x_{i})}$ be effectively the same for given $x_{i}$ and
this does not imply that the unconditional models $M_{i,T_{i}}$ are
effectively the same so we cannot conclude that $SC\subset S.$ We do have,
however, that the conditional inference bases are equivalent under
$S.$\smallskip

\noindent\textbf{Proposition 5. }If $(M_{1},x_{1})$ and $(M_{2},x_{2})$ are
equivalent under $SC,$ then the conditional inference bases $(M_{1,T_{1}\mid
L_{1}(x_{1})},T_{1}(x_{1}))$ and $(M_{2,h(T_{2})\mid L_{2}(x_{2})}%
,h(T_{2}(x_{2})))$ are equivalent under $S.$\smallskip

The following result demonstrates that the evidence function $Ev_{SC}$ is the
ultimate presentation of the evidence based upon $S$ and $SC.\smallskip$

\noindent\textbf{Proposition 6.} For data $x$ and model $M=(\mathcal{X}%
,\{P_{\theta,X}:\theta\in\Theta\}),$ the evidence function defined by
(\ref{stabev}) satisfies $Ev_{SC}=Ev_{MS}\circ Ev_{SC}=Ev_{SC}\circ
Ev_{MS}.\smallskip$

\noindent So, the evidence function that results from the two principles, can
be unambiguously defined as the inference base containing both the observed
value of the mss and the collection of conditional distributions given the
laminal ancillary function of the mss as indexed by the model parameter.

The consequence of this development is that the application of the two
principles can be thought of unambiguously as a function on the set of all
inference bases. It is not clear that there shouldn't be further reductions in
$(\mathcal{T},\{P_{\theta,T|L(x)}:\theta\in\Theta\},T(x))$ to remove
ingredients that are still extraneous to the expression of the evidence
concerning $\theta_{true},$ but at this point it is not obvious what form
those would take.

Also, statistical evidence is ultimately expressed as part of answering
statistical questions. For example, what is the appropriate estimate of
$\psi_{true}=\Psi(\theta_{true})$ and how accurate is it or is there evidence
for or against a hypothesis $H_{0}:$ $\Psi(\theta_{true})=\psi_{0}$ and how
strong is this evidence? Simply stating an inference base does not answer such
questions but at least it does tell us what to focus on when devising the answer.

\section{Conclusions}

Various ambiguities have raised doubts about the possibility of a successful
theory for frequentist inference. For example, Birnbaum's theorem concerning
$S$ and $C$ seemingly implying $L$ or for that matter $C$ alone implying $L$
are but two examples. While the validity of these conclusions has been
challenged, consideration of these results still raises concerns as to what
the correct applications of the principles are. For $S$ this is undoubtedly
discarding all aspects of the inference base that are extraneous to expressing
the evidence about $\theta_{true}$ and this leads to the principle as
expressed by Durbin (1970) together with the evidence function $Ev_{MS}$ which
we add to the development. For $C\,$our thesis is that the fundamental idea
underlying the principle is better expressed by $SC$ and the evidence function
$Ev_{SC}$ as this removes the ambiguity about which ancillary to condition on
and avoids any contradictions in the justification for the irrelevance of the
distribution of the ancillary. While the laminal ancillary may often be
trivial, namely, a function constant on the sample space, it seems clear that
we have to accept the verdict that conditioning on any ancillary other than
the laminal is not appropriate. The results developed here have shown that the
principles $S$ and $SC$ are mutually compatible and satisfy the basic
requirement of any statistical principle by inducing equivalence relations on
the set of all inference bases. As such both the logical and statistical
inconsistences in the definition of $C$ have been avoided.

It is true that the stable conditionality principle proposed here, is - in
part - mathematically supported by the taxonomy results in Basu (1959). The
present paper shows, however, that conditioning on stable ancillaries removes
the logical inconsistencies of the standard conditionality principle and
provides a coherent framework for the assessment of statistical evidence.

Certainly this is not the end of the story concerning the concept of
statistical evidence and how it should be measured and expressed, but our hope
is that clarifying the roles of two key principles contributes to a more solid
foundation for statistics.

\section{Appendix}

\subparagraph{\noindent Proof of Proposition 1}

Suppose $U$ is a strong ancillary for $M$ and let $(p_{1},p_{2},\ldots)$ be an
alternative probability distribution on $%
\mathbb{N}
$ for the marginal distribution of $V.$ Then, summing over those $i$ for which
$P_{V}(\{i\})>0$ (otherwise $P_{\theta,U}(\cdot\,|\,V=i)$ is not defined),%
\begin{align*}
&  \sum_{i:P_{V}(\{i\})>0}p_{i}P_{\theta,U}(B\,|\,V=i)=\sum_{i:P_{V}%
(\{i\})>0}\frac{p_{i}}{P_{V}(\{i\})}P_{\theta,X}(U^{-1}B\cap V^{-1}\{i\})\\
&  =\sum_{i:P_{V}(\{i\})>0}\frac{p_{i}}{P_{V}(\{i\})}\sum_{j\in B}P_{\theta
,X}(V^{-1}\{i\}\,|\,U=j)P_{U}(\{j\})\\
&  =\sum_{i:P_{V}(\{i\})>0}\frac{p_{i}}{P_{V}(\{i\})}\sum_{j\in B}%
P_{V}(\{i\}\,|\,U=j)P_{U}(\{j\})
\end{align*}
where the last equality follows because $U$ is strong which implies that $V$
is ancillary when the mixture distribution for $U$ puts all its mass at $j$ so
$P_{\theta,X}(V^{-1}\{i\}\,|\,U=j)$ is independent of $\theta$ as is the sum.
Therefore, $U$ is stable.

Now suppose $U$ is a stable ancillary and $V$ is ancillary and let
$(p_{1},p_{2},\ldots)$ be an alternative probability distribution on $%
\mathbb{N}
$ for the marginal distribution of $U.$ Then, summing over those $i$ for which
$P_{U}(\{i\})>0,$%
\begin{align*}
&  \sum_{i:P_{U}(\{i\})>0}p_{i}P_{\theta,V}(B\,|\,U=i)=\sum_{i:P_{U}%
(\{i\})>0}\frac{p_{i}}{P_{U}(\{i\})}P_{\theta,X}(V^{-1}B\cap U^{-1}\{i\})\\
&  =\sum_{i:P_{U}(\{i\})>0}\frac{p_{i}}{P_{U}(\{i\})}\sum_{j\in B}P_{\theta
}(U^{-1}\{i\}\,|\,V=j)P_{V}(\{j\})\\
&  =\sum_{i}\frac{p_{i}}{P_{V}(\{i\})}\sum_{j\in B}P_{U}(\{i\}\,|\,V=j)P_{V}%
(\{j\})
\end{align*}
where the last equality follows because $U$ is stable. Therefore, $U$ is
strong. $\blacksquare$

\subparagraph{\noindent Proof of Proposition 2}

(i) Suppose $U$ is an ancillary and it is not a function of a maximal
ancillary $V.$ Then it cannot be that $(U,V)$ is ancillary because, if it
were, then $V$ is a function of $(U,V)$ and thus is not maximal. Since
$P_{\theta,X}(U\in A,V\in B)=P_{\theta,X}(V\in B\,|\,U\in A)P_{U}(A),$ it
cannot be the case that $P_{\theta,X}(V\in B\,|\,U\in A)$ is independent of
$\theta$ for every $A$ and $B$ and so $U$ is not strong. Therefore, any strong
ancillary is a function of every maximal ancillary.

Conversely, suppose $U$ is a minimal ancillary and $V$ is another ancillary.
Then $V$ can be expressed as function of some maximal ancillary $W,$ say
$V=h(W),$ and since $U$ is minimal, it can also be expressed as $k(W)$ for
some function $k,$ Then $P_{\theta,X}(U\in A,V\in B)=P_{\theta,X}(W\in
k^{-1}A\cap h^{-1}B)$ which is independent of $\theta$ because $W$ is
ancillary. Therefore,%
\[
\sum_{i:P_{V}(\{i\})>0}p_{i}P_{\theta,U}(A\,|\,V=i)=\sum_{i:P_{V}%
(\{i\})>0}\frac{p_{i}}{P_{V}(\{i\})}P_{\theta,X}(k^{-1}A\cap h^{-1}\{i\})
\]
which is independent of $\theta$ for every probability distribution
$(p_{1},p_{2},\ldots)$ $%
\mathbb{N}
.$ Therefore, $U$ is a stable ancillary.

\noindent(ii) Let $A_{1},A_{2},\ldots$ be a list of the minimal ancillaries
for model $M$ and put $A=(A_{1},A_{2},\ldots):\mathcal{X}\rightarrow
\mathcal{A}_{1}\mathcal{\times A}_{2}\mathcal{\times}\ldots.$ Now let
$C_{i}\in\mathcal{C}_{i}$ and then $A^{-1}(C_{1}\times C_{2}\times
\ldots)=A^{-1}C_{1}\cap A^{-1}C_{2}\cap\mathcal{\cdots\in B}$ so we can write
$A:(\mathcal{X},\mathcal{B)\rightarrow(A}_{1}\mathcal{\times A}_{2}%
\mathcal{\times}\ldots,\mathcal{C}_{1}\mathcal{\times C}_{2}\mathcal{\times
}\ldots)$ and $A$ is a valid statistic. Further, for a maximal ancillary $W$
there exist functions $h_{1},h_{2},\ldots$ such that $A_{i}=h_{i}(W)$ and this
implies that $A$ is ancillary. For any other ancillary $U,$ there exist a
maximal ancillary $W$ and function $h$ such that $U=h(W)$ and also there are
functions $h_{1},h_{2},\ldots$ such that $A_{i}=h_{i}(W)$ and this implies
that $(A,U)$ is ancillary. As such this proves that $A$ is a minimal ancillary
and moreover it is maximal in this class because every other minimal ancillary
is a function of $A.$ $\blacksquare$

\subparagraph{\noindent Proof of Proposition 3}

We need to show that the relation given by $SC$ is (i) reflexive, (ii)
symmetric and (iii) transitive.

(i) Suppose model $M$ has mss $T$ and laminal $L\in\mathcal{A}_{M}.$ Then
taking $M_{i}=M,T_{i}=T,L_{i}=L$ for $i=1,2$ and $h$ equal to the identity in
(\ref{screl}) establishes reflexivity.

(ii) Symmetry also follows because (\ref{screl}) implies%
\[
(M_{2,T_{2}\mid L_{2}(x_{2})},T_{2}(x_{2}))=(M_{1,h^{-1}(T_{1})\mid
L_{1}(x_{1})},h^{-1}(T_{1}(x_{1}))).
\]

(iii) Finally suppose that $(M_{1},x_{1})$ and $(M_{2},x_{2})$ are related
under $SC$ as well as $(M_{2},x_{2})$ and $(M_{3},x_{3}).$ Let $T_{i},L_{i}%
\in\in\mathcal{A}_{M_{i}}$ denote the mss and laminal ancillaries for $M_{i}$
and $h_{12}:\mathcal{T}_{2}\mathcal{\rightarrow T}_{1},h_{23}:\mathcal{T}%
_{3}\mathcal{\rightarrow T}_{2}$ be the 1-1, onto mappings that are used in
(\ref{screl}) to establish these relations. Then
\begin{align*}
(M_{1,T_{1}\mid L_{1}(x_{1})},T_{1}(x_{1}))  &  =(M_{2,h_{12}(T_{2})\mid
L_{2}(x_{2})},h_{12}(T_{2}(x_{2}))),\\
(M_{1,T_{2}\mid L_{2}(x_{2})},T_{2}(x_{2}))  &  =(M_{3,h_{23}(T_{3})\mid
L_{3}(x_{3})},h_{23}(T_{3}(x_{3})))
\end{align*}
both hold. Now define $h_{13}=h_{12}\circ h_{23}.$ Then if follows that
\begin{align*}
(M_{3,h_{12}\circ h_{23}(T_{3})\mid L_{3}(x_{3})},h_{12}\circ h_{23}%
(T_{3}(x_{3})))  &  =(M_{2,h_{12}(T_{2})\mid L_{2}(x_{2})},h_{12}(T_{2}%
(x_{2})))\\
&  =(M_{1,T_{1}\mid L_{1}(x_{1})},T_{1}(x_{1}))
\end{align*}
and this establishes that $(M_{1},x_{1})$ and $(M_{3},x_{3})$ are related
under $SC$ so the relation is transitive. $\blacksquare$

\subparagraph{\noindent Proof of Proposition 4}

Suppose that $(M_{1},x_{1})$ and $(M_{2},x_{2})$ are equivalent under $S$ so%
\[
(M_{1,T_{1}},T_{1}(x_{1}))=(M_{2,h(T_{2})},h(T_{2}(x_{2}))).
\]
Since the models $M_{1,T_{1}}$ and $M_{2,T_{2}}$ are relabellings of each
other via $h,$ this implies that the ancillarity structure of the two models
is effectively (via the relabelling) the same and, in particular, the laminals
$L_{1}\in\mathcal{A}_{M_{1}}$ and $L_{2}\in\mathcal{A}_{M_{2}}$ are related
via $L_{1}=h(L_{2}).$ This implies%
\[
(M_{1,T_{1}\mid L_{1}(x_{1})},T_{1}(x_{1}))=(M_{2,h(T_{2})\mid L_{2}(x_{2}%
)},h(T_{2}(x_{2})))
\]
and so $(M_{1},x_{1})$ and $(M_{2},x_{2})$ are equivalent under $SC.$
$\blacksquare$

\subparagraph{\noindent Proof of Proposition 5}

Since the two conditional models are simply relabellings it must be the case
that they have effectively the same minimal sufficient statistics and this
implies the result. $\blacksquare$

\subparagraph{\noindent Proof of Proposition 6}

Since $Ev_{SC}(M,x)$ only depends on the model and data through the model for
a mss $T$ and the observed value of $T(x),$ and we have restricted to
ancillaries that are functions of the mss, it is clear that $Ev_{SC}\circ
Ev_{MS}(M,x)=Ev_{SC}(M,x).$

Now consider the reverse order where $Ev_{SC}$ outputs $(\mathcal{T}%
,\{P_{\theta,|L(x)}:\theta\in\Theta\},T(x))$ based on laminal ancillary
$L\in\mathcal{A}_{M}$. We can write $L$ as $L=g(T(x))$ for some function
$g.$The sample space for $T$ in this conditional model is $\{t\in
\mathcal{T}:g(t)=g(T(x))\}$ and which is a union of preimage contours of $T.$
For a $t$ satisfying $g(t)=g(T(x))$ then $P_{\theta,T|L(x)}(\{t\})=P_{\theta
,T}(\{t\})/P_{L}(\{L(x)\}).$ Therefore, if $t_{1},t_{2}$ are distinct elements
of $\{t:g(t)=g(T(x))\},$ then we cannot have
\[
P_{\theta,T|L(x)}(\{t_{1}\})=cP_{\theta,T|L(x)}(\{t_{2}\})
\]
for every $\theta$ for some constant $c>0,$ otherwise we would have
$P_{\theta,T}(\{t_{1}\})=cP_{\theta,T}(\{t_{2}\})$ for every $\theta$ and then
$T$ would not be a mss for the original model. This also implies that the
identity function is a mss for the conditional model which implies
\[
Ev_{MS}(\mathcal{T},\{P_{\theta,T|L(x)}:\theta\in\Theta\},T(x))=(\mathcal{T}%
,\{P_{\theta,T|L(T(x))}:\theta\in\Theta\},T(x)),
\]
namely, there is no reduction. This proves the result. $\blacksquare$

\section{References}

\noindent Basu, D. (1959) The family of ancillary statistics. Sankhy\={a} 21,
247-256.\smallskip

\noindent Basu, D. (1964). Recovery of ancillary information. Sankhy\={a} 26,
3-16.\smallskip

\noindent Birnbaum, A. (1962) On the foundations of statistical inference
(with discussion). \textit{J. Amer. Stat. Assoc.}, 57, 269-332.\smallskip

\noindent Buehler, R. J. (1982). Some ancillary statistics and their
properties. J. Amer. Statist. Assoc. 77, 581-594.\smallskip

\noindent Cox, D. R. (1958) Some problems connected with statistical
inference. The Annals of Mathematical Statistics, 29 (2), 357-372.\smallskip

\noindent Cox, D. R. (1971). The choice between alternative ancillary
statistics. J. Roy. Statist. Soc., B 33, 251-252.\smallskip

\noindent Durbin, J. (1970) On Birnbaum's theorem on the relation between
sufficiency, conditionality and likelihood. \textit{J. Amer. Stat. Assoc.},
654, 395-398.\smallskip

\noindent Evans, M., Fraser, D.A.S. and Monette, G. (1986) On principles and
arguments to likelihood (with discussion). \textit{Canad. J. of Statistics},
14, 3, 181-199.\smallskip

\noindent Evans, M. (2013) What does the proof of Birnbaum's theorem prove?
Electronic Journal of Statistics, Volume 7, 2645-2655.\smallskip

\noindent Evans, M. (2015) Measuring Statistical Evidence Using Relative
Belief. Monographs on Statistics and Applied Probability 144, CRC
Press.\smallskip

\noindent Ghosh M., Reid, N. and Fraser, D. A. S. (2010) Ancillary statistics:
a review. Statistica Sinica 20, 1309-1332.\smallskip

\noindent Kalbfleisch, J.D. (1975) Sufficiency and conditionality.
\textit{Biometrika}, 62, 251-259.\smallskip

\noindent Stigler, S. (2001) Ancillary History. IMS\ Lecture Notes-Monograph
Series, 36, State of the Art in Probability and Statistics, 555-567.

\end{document}